\documentclass[11pt]{revtex4}

\usepackage{graphicx,epsfig,amssymb,amsbsy,amsmath}

\newcommand{\arccot}{\operatorname{arccot}}

\newcommand{\cref}[1]{chapter~\ref{#1}}

\begin{document}

\title{Non-deterministic dynamics of a mechanical system}
\author{Robert Szalai \& Mike R. Jeffrey}\address{Engineering Mathematics, University of Bristol, UK, email: r.szalai@bristol.ac.uk}
\date{\today}

\begin{abstract}
A mechanical system is presented exhibiting a non-deterministic singularity,
that is, a point in an otherwise deterministic system where forward
time trajectories become non-unique. A Coulomb friction force applies
linear and angular forces to a wheel mounted on a turntable. In certain
configurations the friction force is not uniquely determined. When
the dynamics evolves past the singularity and the mechanism slips,
the future state becomes uncertain up to a set of possible values.
For certain parameters the system repeatedly returns to the singularity,
giving recurrent yet unpredictable behaviour that constitutes non-deterministic
chaotic dynamics. The robustness of the phenomenon is such that we
expect it to persist with more sophisticated friction models, manifesting
as extreme sensitivity to initial conditions, and complex global dynamics
attributable to a local loss of determinism in the limit of discontinuous
friction.
\end{abstract}

\maketitle

\section{Introduction}

Classical mechanics accounts for motion of both rigid and deformable
bodies, describing macroscopic processes, often with great accuracy. Even for simple bodies, however, the contact forces between them can be particularly complex, and one must
resort either to intricate fine scale models or crude empirical approximations.
This is the case for both friction and impact. In impact, reaction
forces vary by several magnitudes within the very short time interval
that contact occurs. Models are generally empirical, taking the form
of kinematic constraints given between infinitesimally close times
where contact is made and lost \cite{b00,stewart2000,stronge04}. Friction, on the other
hand, describes a resistance force between surfaces in sustained contact.
The evolution of the friction force over significant time intervals
stands as a singularly enduring problem in classical rigid body mechanics,
as a brief inspection of the extensive literature
reveals. Common empirical force laws still largely resemble those set down by Da Vinci, Amontons and Coulomb (see e.g. \cite{b00,dowson,olsson98,wojewoda08}), of a constant force opposing the direction of motion. This simplicity is in contrast with experimental images that reveal a complex evolution of the contact surface \cite{akay99,bhushan98,bot09,krim,persson}.

The difficulties in obtaining a clear theory of friction are further complicated by the coefficient of friction varying with scale, material, motion, and physical conditions such as temperature and lubrication. Theorists and experimentalists
continue to seek improved contact models by allowing, for example,
impact compliance \cite{champ11painleve,wiercigroch10},
frictional memory \cite{hess91,dawes11}, speed dependence (Stribeck
effect) \cite{wojewoda08,krim}, and microscopic structure such as
asperity deformation at the frictional interface \cite{akay99,albender04}.

Regardless of the chosen model, mechanical laws are reasonably associated
with deterministic dynamics. Reliable models are vital in the study
of tyre deformation in the automotive and aeronautic industries, of
drill strings, turbine rotors, and so on. To this end, all
of the contact models discussed above tend to retain one essential feature of friction and impact, and
that is discontinuity, either from switching between free-flight and
in-contact impact laws, or switching between ``left-motion'' and
``right-motion'' friction laws. It is easy to assume that one may
smooth out the discontinuity. The effects of smoothing on the dynamics, however, remain
poorly understood, with both theory and simulations suggesting more intricate dynamics organised, nevertheless, around the quantitative features of the discontinuous limit \cite{j13error}.

The crucial point is that existence and
uniqueness of the dynamical solutions of differential equations is not guaranteed when they suffer discontinuities.
In many situations, including all those listed above, the existence
of solutions can be guaranteed at the possible expense of loss of
uniqueness. The most obvious example is a slipping block that comes
to rest on a surface. It is easy to predict where the block will stop, but when it comes to rest information about its previous motion is lost to the environment, and piecing together its previous motion is no longer possible.
This loss of unique reversibility is common in a piecewise-smooth dynamical system. When a
trajectory encounters a discontinuity there are local rules, physically
motivated, that provide a well-defined and solvable system. But those
solutions may be non-unique either in forward or backward time. Backward
time ambiguity corresponds to the onset of sticking that brings a
block to rest, consistent with centuries of physical intuition. Forward
time ambiguity, on the other hand, signals a breakdown in the deterministic
character of physical law.

A classic example of a forward time ambiguity is the Painlev\'e paradox (see e.g.
\cite{champ11painleve}) responsible for the juddering of a piece of
chalk when pushed, instead of dragged, along a blackboard. As a rod
impacts a surface its endpoint can stick in place or slip. The effects
of impact and friction, together with the coupling of the rod's linear
and angular motion, may combine such that the surface seems to attract,
rather than resist, the rod. The rod's release is then not uniquely
determined by the governing equations, resulting in unpredictable
skipping of its endpoint along the surface. In this paper we present
a phenomenon that is potentially more common, and more extreme, where
forward and backward time ambiguity combine to generate chaotic dynamics.

This behavior was described in a theoretical context in \cite{cj10},
and in an abstract model proposed in \cite{j11prl}. The discontinuity
in a piecewise-smooth system can create a so-called two-fold singularity,
where a system switches between two states, which both evolve instantaneously
along the switching locus. Two-fold singularities arise in generic
systems of more than two dimensions (with simpler degenerate forms
in planar systems). They have different forms, some of which are benign,
either in the sense that they are invisible to the dynamics, or that
forward time ambiguities are short-lived because trajectories quickly
recombine. In another form, trajectories enter the singularity
both forwards and backwards in time, traveling (counterintuitively)
from attracting regions to repelling regions of phase space, and in a finite amount of time. The result
is a set of allowed future and past trajectories, with none more likely
than any other. In \cite{cj10},
local conditions near a bifurcation were shown to re-inject trajectories
so that this unpredictable singularity was visited repeatedly, giving
birth to a non-deterministic form of chaos. In \cite{j11prl},
a global re-injection mechanism was shown to create more robust non-deterministic
chaos in an abstract system with negative damping. Attempts have also
been made to characterize the possible appearance of the singularity
in electrical control models \cite{cj09,col11}.

Despite these steps, no testable examples of this singular
prediction of piecewise-smooth dynamical systems theory have previously been devised. In this paper the singularity
is shown to occur in a simple and realistic mechanical model, motivated
from devices such as \cite{TakacsHogan}, involving a Coulomb friction law and
a simplified version of Pacejka's magic formula \cite{pacejka06}
to characterize nonlinear response of a wheel. Non-deterministic chaos
ensues for a range of parameters, and while a continuum model such
as the stretched string tyre model \cite{stepan09} should resolve
the singularity, the sensitivity to initial conditions that gives
rise to chaos-like dynamics is expected to persist.

Below we present the mechanical model in section \ref{sec:model}, and show how a singularity arises in the friction force in section \ref{sec:fsing}. Some generalities of the singularity are presented briefly in section \ref{sec:teixeira}, followed by identification of parameter regimes for which non-deterministic chaos can appear along with simulations in section \ref{sec:sim}. We remark on special simplifying cases in \ref{sec:special} and make concluding remarks in section \ref{sec:conc}.

\section{Mechanical model}\label{sec:model}

The crucial component of the following mechanical system is the friction torque and the friction force imparting on a wheel that can balance each other with respect to the contraints built into the system. The precise design of the mounting of the wheel on a rotating disc is contrived solely to simplify the analysis by eliminating the rotation angle of the mount from the equations of motion.

\begin{figure}[h!]
\begin{centering}
\includegraphics[width=0.8\textwidth]{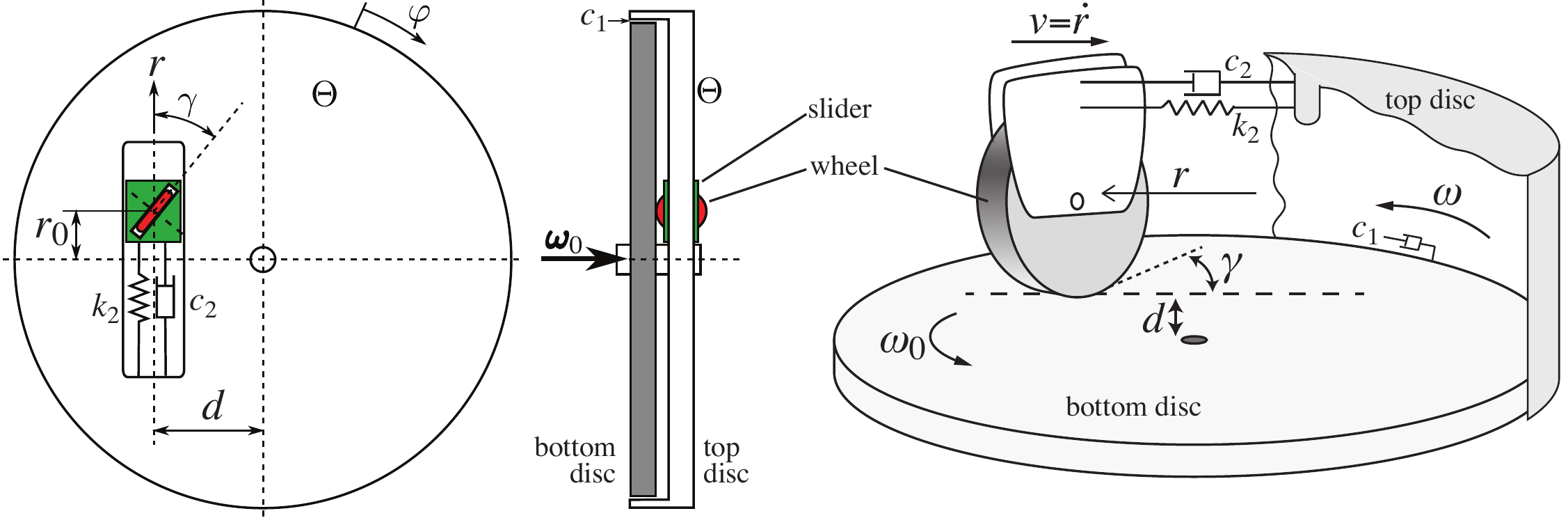}\vspace{-0.11cm}
\caption{\sf Mechanical model of the system show from the top (left)
and from the side (middle), with a conceptual sketch (right). The bottom disc is rotating with
a constant angular speed $\omega_{0}$. The top disc rotates around
the same centre pin, but it is connected to the bottom disc through
a viscous friction coefficient $c_{1}$. A slider is fitted in the top disc, forming a single degree of freedom
oscillator. The slider contains a wheel mounted at an angle $\gamma$ to the
slider direction. The wheel is in contact with the bottom disc and Coulomb friction acts between them.}\label{fig:mechmod}
\end{centering}
\end{figure}

Figure \ref{fig:mechmod} shows the mechanical model of interest.
It consists of two discs that can independently rotate on a shared
axis, and a wheel (shown red) mounted transversely inside the upper
disc so that it rolls over the surface of the lower disc. The wheel
is mounted inside the top disc by means of a slider (a one degree
of freedom oscillator, shown in green) a distance $d$ from the disc
axis. The wheel sits at an angle $\gamma$ to the slider direction.
The slider is fixed to the disc via a spring coefficient $k_{2}$
and damping coefficient $c_{2}$, and it's equilibrium position lies
a distance $\sqrt{d^{2}+r_{0}^{2}}$ from the disc axis. The wheel's
mass is negligible, the slider's mass is $m$, and the moment of inertia
of the top disc is $\Theta=\beta^{2}m$. The bottom disc rotates with
a constant angular speed $\omega_{0}$, with a viscous friction coefficient
$c_{1}$ between it and the upper disc. All of the parameters so far
described are constants.

The most crucial feature of the model is the friction force between
the wheel and the steadily rotating bottom disc. The friction force
acts perpendicular to the plane of the wheel and has magnitude
$\mu$. Inhomogeneities in the force distribution on the contact patch
between the bottom disc and the wheel create a torque on the wheel
of magnitude $\mu M$. The function $M$ models the deformation of
a tyre around the wheel, for which different models are available
\cite{pacejka06,Kuiper07,stepan09}. Most of these tyre models are
smooth but include steep gradients
that we replace with discontinuities. In this paper we choose the
stretched-string tyre model described in Appendix B.

This system can be described using two generalized coordinates, namely the
rotation angle of the top disc, $\varphi$, and the displacement of
the slider, $r$. Differential equations for the slider displacement
$r$, slider speed $v=\dot{r}$, and the top disc's angular speed
$\omega=\dot{\varphi}$ , can be derived from the Euler-Lagrange equations
(or by directly calculating linear and angular accelerations of the
slider and the top disc respectively). These are found to be independent
of the angle variable $\varphi$, giving the three dimensional dynamical
system
\begin{eqnarray}
\dot{r} & = & m\left(\beta^{2}+r^{2}\right)v,\label{ds1}\\
\dot{v} & = & \left(\beta^{2}+d^{2}+r^{2}\right)p_{1}(r,v,\omega)-d\left(c_{1}+2mrv\right)\omega+F(h,g)p_{2}(r)+d\, M(h,g)\label{ds2}\\
\dot{\omega} & = & d\, p_{1}(r,v,\omega)-\left(c_{1}+2mrv\right)\omega+F(h,g)r\cos\gamma+M(h,g),\label{ds3}
\end{eqnarray}
in terms of functions $F,M,p_1,p_2,h,g$ defined below, and where the dot denotes the time derivative. In the equations of motion (\ref{ds1})-(\ref{ds3}), time has been rescaled so that the strictly positive quantity $m\left(\beta^{2}+r^{2}\right)$
does not appear in the denominator of the righthand side. The functions $p_{1}$
and $p_{2}$ are given by
\[
p_{1}(r,v,\omega)=k_{2}(r_{0}-r)-c_{2}v+mr\omega^{2},\; \qquad
p_{2}(r)=dr\cos\gamma+\left(\beta^{2}+r^{2}\right)\sin\gamma.
\]

The functions $g$ and $h$ are the components of the wheel's relative velocity with respect to the bottom disc, in the
rolling and the lateral directions, respectively. They are given by
\begin{align}
h(r,v,\omega) & =-(v-d(\omega-\omega_{0}))\sin(\gamma)-r(\omega-\omega_{0})\cos(\gamma),\label{h}\\
g(r,v,\omega) & =-(v-d(\omega-\omega_{0}))\cos(\gamma)+r(\omega-\omega_{0})\sin(\gamma).\label{g}
\end{align}
When $h=0$ and $v=0$,
the wheel rolls around the disc in a circle without slipping, such
that $r=d\tan\gamma$. When $h=0$ but $v\neq0$, the wheel instantaneously
rolls on a spiral about the centre of the disc.

The friction
force $F$ and moment $M$ represent the tyre model (derived in Appendix B), and are given by
\begin{align}
F(h,g) & =\mathrm{sign}(h)\mu\left(\frac{4}{3}-\frac{\cot\psi}{18\kappa}\right),\label{eq:force}\\
M(h,g) & =\mathrm{sign}(g)\mathrm{sign}(h)\frac{\mu}{6}\cot\psi,\label{eq:moment}
\end{align}
where $\psi=\arccot6\kappa+\left|1-\frac{2}{\pi}\arccot6\kappa\right|\arctan\left|{h}/{g}\right|$. An important feature of the model is that
$\lim_{h\to0\pm}M=\pm\mathrm{sign}(g)\mu\kappa$, which imparts a
nonzero torque on the wheel as it enters the sticking phase $h=0$.

The system is piecewise-smooth because the friction force $F$ and
moment $M$ are discontinuous at $h=0$. This defines a set of points
\[
\Sigma=\{(r,v,\omega)\in\mathbb{R}^{3}:h(r,v,\omega)=0\}
\]
called the \emph{switching surface}. When $h\neq0$, $F$ and
$M$ are smooth functions of $h$ and $g$, so equations (\ref{ds1})-(\ref{ds3})
remain smooth and have unique solutions $\left(r(t),v(t),\omega(t)\right)$.
Note that the moment $M$ stays smooth when $g$ changes sign since this occurs when
$\psi=\pi/2$ and hence $\cot\psi=0$.

To solve the equations of motion when they reach the switching surface $\Sigma$, we first assume that $F$ and $M$ never exceed the absolute values they have as they approach $\Sigma$ (i.e. static and kinetic friction are equal), hence $|F|\le\mu$ and $|M|\le\mu\kappa$. Only two kinds of solutions are possible given that a trajectory $(r(t),v(t),\omega(t))$ must be a continuous solution of (\ref{ds1})-(\ref{ds2}): either a solution will stick to the surface $\Sigma$, or it will  cross through $\Sigma$ transversally which reverses the direction of slipping.

The values of $F$ and $M$ during sticking are not fixed by (\ref{eq:force})-(\ref{eq:moment}) because the value of ${\rm sign}(h)$ is not well-defined at $h=0$. To establish whether sticking occurs, we ask whether there exists $F\in[-\mu,+\mu]$ such that the vector field (\ref{ds1})-(\ref{ds3}) is tangential to $\Sigma$. If such an $F$ exists, then sticking occurs and solution follow (\ref{ds1})-(\ref{ds3}) along $\Sigma$ with the given value of $F$ and with $M=\kappa F$. In mechanical terms, evolution along $\Sigma$
corresponds to stick of the contact surfaces, during which the mechanical constraint $h=0$ holds. Note that the friction force and moment are
tied together during stick by the constant $\kappa$.
If there exists no value of $\lambda$ in the interval $[-\mu,+\mu]$ for which (\ref{ds1})-(\ref{ds3}) is tangential to $\Sigma$, then sticking to $\Sigma$ is impossible dynamically and trajectories cross from slipping in one direction to the other.

\begin{figure}[h!]
\begin{centering}
\includegraphics[width=0.7\textwidth]{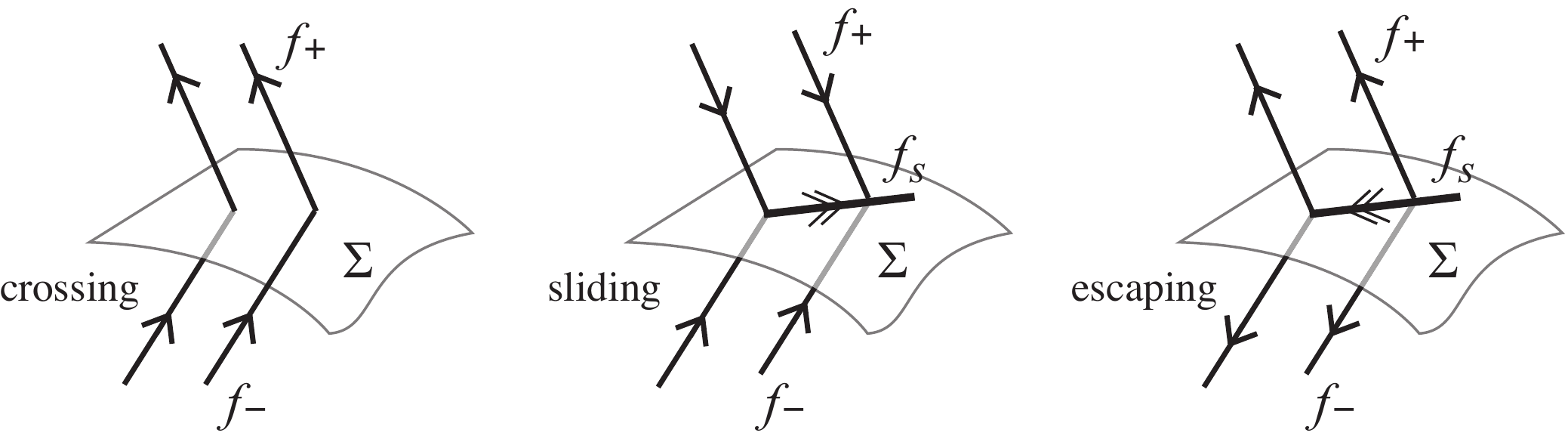}\vspace{-0.11cm}
\caption{\sf Typical dynamics as the switching surface: crossing, sliding, and escaping. The limit of (\ref{ds1})-(\ref{ds2}) as $(r,v,\omega)$ approaches $\Sigma$ are labelled as $f_\pm=\lim_{h\rightarrow0\pm}(\dot r,\dot v,\dot\omega)$. The vector field $f_s$ that governs sticking will be defined in section \ref{sec:teixeira}.}\label{fig:fil}\end{centering}
\end{figure}

The region of
$\Sigma$ where sticking occurs is reachable from both `right slip' $h>0$ and `left slip' $h<0$. It is called the \emph{sliding surface}, because trajectories do not just stick to $\Sigma$ but evolve or {\it slide} over its surface as $F$ varies between $[-\mu,+\mu]$ (to avoid confusion note: we employ this mathematical usage of the word `slide' during the stick phase on $h=0$, and use `slip' to describe non-stick contact in the mechanical sense for $h\neq0$). The sticking surface a subset of $\Sigma$ on which the vector field $(\dot x,\dot v,\dot\omega)$ points towards $h=0$ from both $h>0$ and $h<0$., hence it is an attractor of the local dynamics. Conversely, it is mathematically possible for an unstable form of sticking to occur on an \emph{escaping surface}, where the vector field $(\dot x,\dot v,\dot\omega)$ points away from $h=0$ from both $h>0$ and $h<0$. Thus an escaping region is a repellor of local dynamics, on which sticking occurs but can end spontaneously (or under arbitrarily small perturbations in practice) upon which solutions depart $\Sigma$ abruptly. This definition seems to imply that solutions
will never reach an escaping surface, making it physically uninteresting. That this it not the case was shown in \cite{cj10}, and we illustrate it shortly using the model above. The three basic behaviours --- sliding, escaping, and crossing --- are illustrated in Figure \ref{fig:fil}.

Sticking ceases at points where the attractivity of $\Sigma$ breaks down, namely where $(\dot r,\dot v,\dot\omega)$ does not point either towards $\Sigma$ (for a sliding region) or away from $\Sigma$ (for an escaping region) from both $h>0$ and $h<0$. Thus the conditions
$$\lim_{h\rightarrow0\pm}(\dot r,\dot v,\dot\omega)\cdot\nabla h(r,v,\omega)=0\;$$
define curves bounding regions of sliding or escaping.
If these two curves intersect they typically bring sliding and escaping regions together at the intersection point, defined by
\begin{equation}\label{2fold}
\lim_{h\rightarrow0+}(\dot r,\dot v,\dot\omega)\cdot\nabla h(r,v,\omega)=\lim_{h\rightarrow0-}(\dot r,\dot v,\dot\omega)\cdot\nabla h(r,v,\omega)=0\;.
\end{equation}
Despite its low dimension, this point, known as the {\it two-fold singularity}
, will turn out to be highly important to the dynamics. One role it can have is to funnel trajectories from sliding to escaping motion \cite{cj10}, making escaping dynamics an important and observable feature of the system's dynamics. This is the situation found in the device above. Once inside the escaping region, it is impossible to determine from the equations of motion exactly how and where
a trajectory will cease sticking. Thus the singularity introduces an interval of non-deterministic motion in an otherwise deterministic system, in which we cannot determine how long sticking will last, or what the slipping motion will look like immediately after, but once slipping motion is restarted then determinism is restored.

These features of crossing and sticking will be formalised in more detail in section \ref{sec:teixeira}.

\section{Singularity of the friction force}\label{sec:fsing}

The singularity that appears is simple to describe in
mechanical terms if we assume that the friction
force and moment are related linearly during stick. If the mechanical
system can be placed into a configuration where friction force and
moment balance each other with respect to the constraints, then due
to their linear relation any friction force would balance the moment
that is proportional to the force. In the model (\ref{ds1})-(\ref{ds3}), this balance occurs
via the constraint of the pin holding together the two
independently rotating discs. The formal equation (\ref{2fold})
exactly describes this situation. It requires stick at a point with
two different values of $\lambda$, which due to linearity implies
that (\ref{2fold}) must also hold for any $\lambda$, hence
$M$ and $F$ are in balance.

The conditions $h=0$ and (\ref{2fold}) for the two-fold singularity give equations that are solveable, but highly nonlinear, in the variables $(r,v,\omega)$. Rather than solving for these in terms of the system parameters, let us say that the singularity lies at coordinates $(r,v,\omega)=(r^*,v^*,\omega^*)$ when two of the parameters of the device, say $\kappa$ and $k_2$, take certain values dependent on $(r^*,v^*,\omega^*)$. We thus solve the equations (\ref{2fold}) along with $h(r^*,v^*,\omega^*)=0$, to find
\begin{align}
v^* & =\left(\omega^*-\omega_{0}\right)\left(d-r^*\cot\gamma\right),\label{eq:stickV}\\
\kappa & =-\frac{2{r^*}^{2}+\beta^{2}\left(1-\cos(2\gamma)\right)}{2r^*\cos\gamma},\label{eq:singKappa}\\
k_{2} & =\frac{1}{\left(r^*-r_{0}\right)p_{2}(r^*,v^*,\omega^*)}\Bigl(\sin(\gamma)\left(\beta^{2}+r^{\star2}\right)\left(mr^*{\omega^*}^2-c_{2}v^*\right)\nonumber \\
 & \quad\quad-\cos(\gamma)\left(c_{1}r^*\omega^*+c_{2}dr^*v^*+m\left({r^*}^2\left(v^*(\omega^*+\omega_{0})-d{\omega^*}^2\right)+\beta^{2}v^*(\omega_{0}-\omega^*)\right)\right)\Bigr).\label{eq:singK2}
\end{align}
By solving these numerically one finds at most four singularities on $\Sigma$. Only one of these creates non-determinism, as described in the next section.

\section{Teixeira singularity}\label{sec:teixeira}

\begin{figure}[h!]
\begin{centering}
\includegraphics[width=0.9\textwidth]{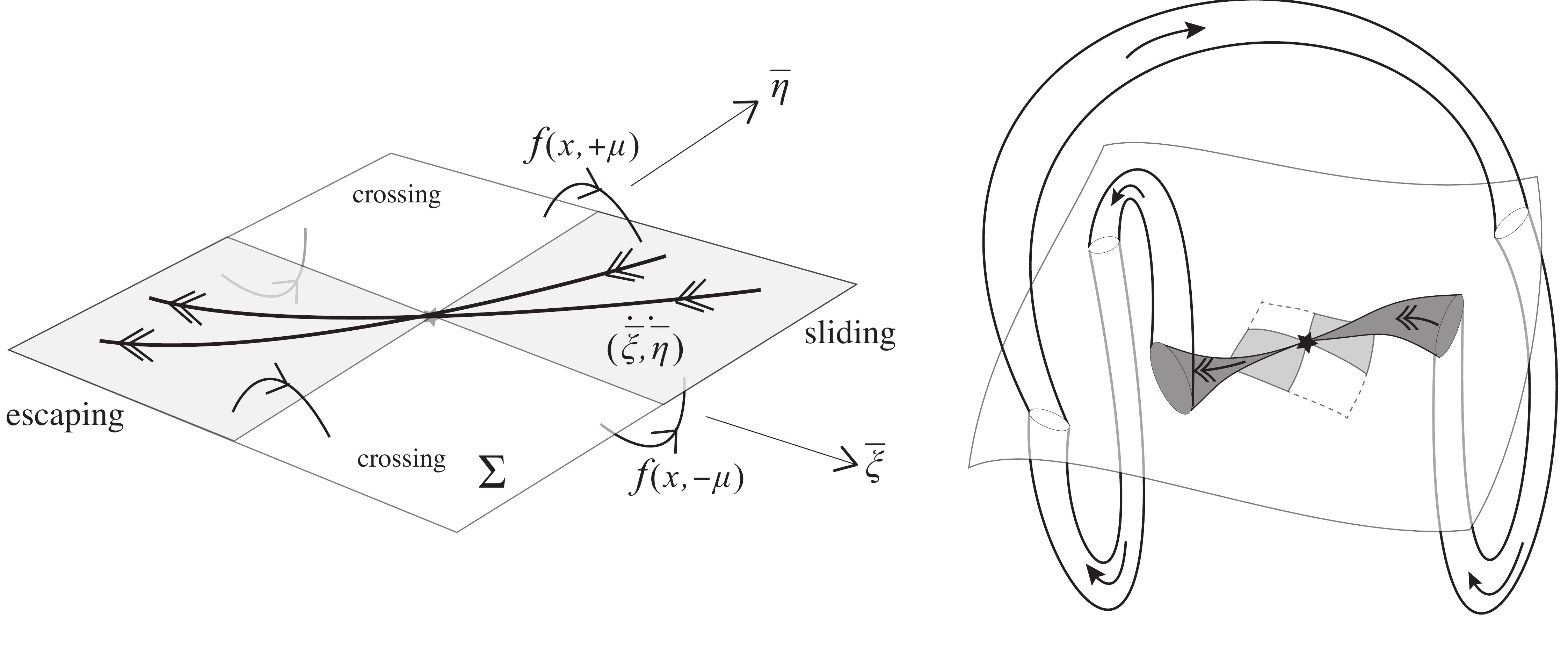}\vspace{-0.11cm}
\caption{\sf Non-determinism at a Teixeira singularity. {\it Left}: Slipping (single-headed arrows) follows $f(x,+\mu)$ or $f(x,-\mu)$ either side of the switching surface. Sticking (double arrows) follows $(\dot{\bar\xi},\dot{\bar\eta})$ in the sliding and escaping regions (shaded). In the case shown, the sticking flow passes through the singularity non-deterministically.
{\it Right}: Sketch of non-deterministic chaos in the mechanical system. The small
shaded ``bow-tie'' shows the local region of the singularity (shown left). A set of trajectories loop around via a sequence
of slipping (white bands) with switches, connecting to the singularity
via sticking (dark band) in both forward and backward time.}\label{fig:2fold}\end{centering}
\end{figure}

There are different types of two-fold singularities \cite{f88,cj10,jh09}. The type responsible for introducing non-determinism to the device above is known as a {\it Teixeira singularity}, formed where the flow curves inwards towards the switching surface $\Sigma$ from both sides (see later in Figure \ref{fig:2fold}(left)); in the two other types of two-folds, the flow curves away from the surface on one or both sides. Here we derive the conditions that guarantee a non-degenerate singularity, and moreover give
stick motion through the singularity, leading to loss of determinism.

The value of the function ${\rm sign}(h)$ in (\ref{eq:force})-(\ref{eq:moment}) is ambiguous at $h=0$. To resolve this we replace the contact force $\mu{\rm sign}(h)$ with a variable $\lambda$ satisfying
\begin{equation}\label{sign}
\lambda(r,v,\omega)\;\;\in\;\;\mu\times\left\{\begin{array}{lll}{\rm sign}(h(r,v,\omega))&\rm if&h(r,v,\omega)\neq0\;,\\\left[-1,+1\right]&\rm if&h(r,v,\omega)=0\;,\end{array}\right.
\end{equation}
and seek conditions to define the value of $\lambda$ during sticking on $h=0$. Firstly, note that on $h=0$ the friction force and moment are
tied together as $F=\lambda$ and $M=\kappa\lambda$.

For conciseness let $x=(r,v,\omega)$, and denote the right-hand side of (\ref{ds1})-(\ref{ds3}) by $f(x,\lambda)$, giving
\begin{equation}
\dot{x}=f(x,\lambda).\label{eq:VF}
\end{equation}
The vector field \eqref{eq:VF} depends linearly
on the friction force $\lambda$, and we can write
\[
\dot{x}=f(x,0)+\lambda f_{\lambda}(x,0),
\]
where $f_{\lambda}(x,0)$ denotes the partial derivative $\frac{\partial\;}{\partial\lambda}f(x,\lambda)$, evaluated at $\lambda=0$. Paraphrasing the statements at the end of section \ref{sec:model}, we can then find the conditions defining sliding or escaping regions, and find the equations of sticking motion within them.

For sticking motion to occur, the vector field $f(x,\lambda)$ must be able to lie tangent to $\Sigma$ for admissible values of $\lambda$, namely $\lambda\in[-\mu,+\mu]$. This means the $h$ component of $f$ must vanish, or equivalently that the Lie derivative of $h$ with respect to $f$, must vanish,
\begin{equation}
h_{x}(x)\cdot\left(f(x,0)+\lambda f_{\lambda}(x,0)\right)=0,\label{eq:slidCond}
\end{equation}
with solution
\begin{equation}
\lambda=-\frac{h_{x}(x)\cdot f(x,0)}{h_{x}(x)\cdot f_{\lambda}(x,0)}.\label{eq:muCalc}
\end{equation}
Using this value for $\lambda$ the equations of motion for sticking become
\begin{equation}
\dot{x}=f(x,0)-\frac{h_{x}(x)\cdot f(x,0)}{h_{x}(x)\cdot f_{\lambda}(x,0)}f_{\lambda}(x,0).\label{eq:svf}
\end{equation}
This is known as the {\it sliding vector field} in piecewise-smooth dynamical systems theory \cite{f88,bc08}.
In mechanical terms, evolution along $\Sigma$
corresponds to stick of the contact surfaces, during which the mechanical constraint $h=0$ holds. (We re-iterate that the word `slide' is used to describe evolution on $h=0$ during the stick phase, and `slip' is used to describe non-stick contact in the mechanical sense for $h\neq0$). Then $\lambda$ is the static friction force, as calculated from the equations of motion.
If there are no values of $\lambda$ in the interval $[-\mu,+\mu]$ such that (\ref{eq:slidCond}) is satisfied, then sticking to $\Sigma$ is impossible dynamically and solutions must cross $\Sigma$ transversally from left slip, $h<0$, to right slip, $h>0$.

A sliding region, where sticking occurs so (\ref{eq:muCalc}) lies inside $[-\mu,+\mu]$ and $\Sigma$ is attractive, now satisfies
\begin{equation}\label{slregion}
-h_{x}(x)\cdot f_{\lambda}(x,0)>h_{x}(x)\cdot f(x,0)/\mu>h_{x}(x)\cdot f_{\lambda}(x,0)\;.
\end{equation}
An escaping region, where sticking occurs and $\Sigma$ is repulsive, satisfies
\begin{equation}\label{esregion}
-h_{x}(x)\cdot f_{\lambda}(x,0)<h_{x}(x)\cdot f(x,0)/\mu<h_{x}(x)\cdot f_{\lambda}(x,0)\;.
\end{equation}
The boundaries of sliding or escaping regions hence lie where
$$h_{x}(x)\cdot f(x,0)=\pm\mu h_{x}(x)\cdot f_{\lambda}(x,0)\;.$$
If both of these equations hold, then two boundaries intersect at a point, which we label $x^*$, where
\begin{equation}\label{2foldcond}
h_{x^*}(x^*)\cdot f(x^*,0)= h_{x^*}(x^*)\cdot f_{\lambda}(x^*,0)=0\;.
\end{equation}
This is the two-fold singularity. Generically the intersection of the two boundaries is transversal, meaning $f(x^*,0)$ and $f_{\lambda}(x^*,0)$
are linearly independent, as well as both tangent to $\Sigma$.
At $x^*$, both the numerator and the denominator of the fraction in (\ref{eq:svf})
vanish, and so the sliding vector field is ill-defined. The problem is easily dealt with by considering the scaled vector field
\begin{equation}
\left(h_{x}(x)\cdot f_{\lambda}(x,0)\right)\dot{x}=h_{x}(x)\cdot f_{\lambda}(x,0)f(x,0)-h_{x}(x)\cdot f(x,0)f_{\lambda}(x,0).\label{eq:slidRescale}
\end{equation}
The right-hand side is a vector field which is a positive scaling of (\ref{eq:svf}) in the sliding region, a scaling of (\ref{eq:svf}) with directions reversed in the escaping region, and with a zero at the singularity (see e.g. \cite{cj10}).
We are interested in the dynamics local to the singularity, so we
linearize equation \eqref{eq:slidRescale} at $x^*$, which yields
\begin{equation}
\left(\partial_{x}\left(h_{x}\cdot f_{\lambda}\right)x\right)\dot{x}=f\;\partial_{x}\left(h_{x}\cdot f_{\lambda}\right)x-f_{\lambda}\partial_{x}\left(h_{x}\cdot f\right)x+O(|x|^2)\label{eq:linSlide}
\end{equation}
(omitting the arguments of $h=h(x)$, $f=f(x,0)$ and $f_\lambda=f_\lambda(x,0)$ for brevity).

The vector field on the right side of \eqref{eq:linSlide}
spans two dimensions (the tangent plane to $\Sigma$), since $f(x^*,0)$ and $f_{\lambda}(x^*,0)$
are linearly independent. Therefore it is sensible to parameterize
the perturbation as $x=\alpha f(x,0)+\beta f_{\lambda}(x,0)+y$, where $y$
can be chosen from the kernel of \eqref{eq:linSlide}, and $\alpha,\beta,$ are parameters. In particular
there is a projection $P$ such that $y=Px$.
Let us define
\begin{equation}
\begin{array}{ll}
\mathcal{K}^{++} =\partial_{x}\left(h_{x}\cdot f(x,+\mu)\right)\cdot f(x,+\mu),&\quad\mathcal{K}^{+-}=\partial_{x}\left(h_{x}\cdot f(x,+\mu)\right)\cdot f(x,-\mu),\\
\mathcal{K}^{-+} =\partial_{x}\left(h_{x}\cdot f(x,-\mu)\right)\cdot f(x,+\mu),&\quad\mathcal{K}^{--}=\partial_{x}\left(h_{x}\cdot f(x,-\mu)\right)\cdot f(x,-\mu),
\end{array}
\end{equation}
The boundaries of the sticking regions are found by substituting $\lambda=\pm\mu$ into (\ref{eq:slidCond}). Linearizing
this constraint we find
\begin{equation}
\partial_{x}\left(h_{x}\cdot\left(f(x,0)\pm\mu f_{\lambda}(x,0)\right)\right)\left(\alpha^{\pm}f(x,0)+\beta^{\pm}f_{\lambda}(x,0)\right)+O(|x|^2)=0.\label{eq:TanLines}
\end{equation}
 To leading order this is satisfied if
\begin{align*}
\alpha^{\pm} & =\;\partial_{x}\left(h_{x}\cdot\left(f(x,0)\pm\mu f_{\lambda}(x,0)\right)\right)f_{\lambda}(x,0)\;\;\;=\;\frac{1}{2\mu}\left(\mathcal{K}^{\pm+}-\mathcal{K}^{\pm-}\right),\\
\beta^{\pm} & =-\partial_{x}\left(h_{x}\cdot\left(f(x,0)\pm\mu^{\pm}f_{\lambda}(x,0)\right)\right)f(x,0)=-\frac{1}{2}\left(\mathcal{K}^{\pm-}+\mathcal{K}^{\pm+}\right).
\end{align*}
The pair of vectors $v^{\pm}=\alpha^{\pm}f(x,0)+\beta^{\pm}f_{\mu}(x,0)$
can now be defined as natural coordinates for the linearised sticking motion, that is, $x=\xi v^{+}+\eta v^{-}+y$ where $\xi$ and $\eta$ are coordinates measured along the $v^+$ and $v^-$ directions. In these coordinates
the linearization of the time-scaled sliding vector field \eqref{eq:linSlide} becomes
\begin{equation}
\left(\xi+\eta\right)\left(\dot{\xi}v^{+}+\dot{\eta}v^{-}+\dot{y}\right)=\left(f(x,0)(\xi+\eta)+\mu f_{\lambda}(x,0)(\xi-\eta)\right).\label{eq:trafSlid}
\end{equation}
We also evaluate the time derivatives to get
\begin{equation}\begin{array}{rcl}
\dot{\xi}v^{+}+\dot{\eta}v^{-}&=&\frac{1}{2\mu}\Bigl\{ \;\;f(x,0)\left(\left(\mathcal{K}^{++}-\mathcal{K}^{+-}\right)\dot{\xi}+\left(\mathcal{K}^{-+}-\mathcal{K}^{--}\right)\dot{\eta}\right)\\
&&\quad-\mu f_{\lambda}(x,0)\left(\left(\mathcal{K}^{+-}+\mathcal{K}^{++}\right)\dot{\xi}+\left(\mathcal{K}^{--}+\mathcal{K}^{-+}\right)\dot{\eta}\right)\Bigr\}.\label{eq:trafDeri}
\end{array}\end{equation}
Substituting \eqref{eq:trafDeri} into \eqref{eq:trafSlid} and solving
for the derivatives $\dot{\xi},\dot{\eta},\dot{y}$ yields
\begin{align}
\dot{\xi} & =-\frac{2\mu}{a(\xi,\eta)}\left(\mathcal{K}^{-+}\xi+\mathcal{K}^{--}\eta\right),\nonumber \\
\dot{\eta} & =\frac{2\mu}{a(\xi,\eta)}\left(\mathcal{K}^{++}\xi+\mathcal{K}^{+-}\eta\right),\label{eq:preNormForm}\\
\dot{y} & =0,\nonumber
\end{align}
where $a(\xi,\eta)=\left(\mathcal{K}^{-+}\mathcal{K}^{+-}-\mathcal{K}^{++}\mathcal{K}^{--}\right)\left(\eta+\xi\right)$.

Equation \eqref{eq:preNormForm} describes the dynamics
on the switching surface. We must also establish how a trajectory might arrive at the switching
surface. Let us assume that the two-fold singularity of interest is a Teixeira singularity \cite{jc09,cj10}, meaning the flows curve towards $\Sigma$ at the boundaries of sticking (figure \ref{fig:2fold}). This implies
\begin{align}
\mathcal{K}^{++}<0<\mathcal{K}^{--}.\label{eq:Invisible}
\end{align}
The sliding region satisfies (\ref{slregion}), whose linearization is $\;\partial_{x}\left(h_{x}\cdot f(x,\mu)\right)x<0<\partial_{x}\left(h_{x}\cdot f(x,-\mu)\right)x$. This means sliding occurs for values of
$\xi$ and $\eta$ that satisfy
\begin{equation}
\frac{\eta}{2\mu}\left(\mathcal{K}^{-+}\mathcal{K}^{+-}-\mathcal{K}^{++}\mathcal{K}^{--}\right)<0<\frac{-\xi}{2\mu}\left(\mathcal{K}^{-+}\mathcal{K}^{+-}-\mathcal{K}^{++}\mathcal{K}^{--}\right).\label{eq:SlidCond}
\end{equation}
Similarly escaping occurs for values
\begin{equation}
\frac{-\xi}{2\mu}\left(\mathcal{K}^{-+}\mathcal{K}^{+-}-\mathcal{K}^{++}\mathcal{K}^{--}\right)<0<\frac{\eta}{2\mu}\left(\mathcal{K}^{-+}\mathcal{K}^{+-}-\mathcal{K}^{++}\mathcal{K}^{--}\right).\label{eq:EscCond}
\end{equation}
Two crucial constants characterize the local dynamics,
\[
\mbox{\ensuremath{\mathcal{J}}}_{1}=\frac{\mathcal{K}^{-+}}{\sqrt{-\mathcal{K}^{++}\mathcal{K}^{--}}},\qquad\text{and}\qquad\mathcal{J}_{2}=-\frac{\mathcal{K}^{+-}}{\sqrt{-\mathcal{K}^{++}\mathcal{K}^{--}}}.
\]
If we choose scaled coordinates $(\bar\xi,\bar\eta)=(\xi,\eta)\sqrt{\frac{-\mathcal{K}^{++}}{\mathcal{K}^{--}}}$,
and rescale time by $\bar t=t\sqrt{-\mathcal{K}^{++}\mathcal{K}^{--}}$, the local system simplifies to
\begin{equation}
\left(\begin{array}{c}
\dot{\bar{\xi}}\\
\dot{\bar{\eta}}
\end{array}\right)=\frac{1}{\left(\mbox{\ensuremath{\mathcal{J}}}_{1}\mbox{\ensuremath{\mathcal{J}}}_{2}-1\right)\left(\bar{\eta}+\bar{\xi}\right)}\left(\begin{array}{cc}
\mbox{\ensuremath{\mathcal{J}}}_{1} & 1\\
1 & \mathcal{J}_{2}
\end{array}\right)\left(\begin{array}{c}
\bar{\xi}\\
\bar{\eta}
\end{array}\right).\label{eq:SlidNormForm}
\end{equation}
Using condition \eqref{eq:SlidCond} we find that if $\mbox{\ensuremath{\mathcal{J}}}_{1}\mbox{\ensuremath{\mathcal{J}}}_{2}-1>0$ then the sliding region is simply $\xi>0,\eta>0$, and $\xi<0,\eta<0$ otherwise. A simple time-rescaling guarantees the former case, transforming  \eqref{eq:SlidNormForm} into
\[
\left(\begin{array}{c}
\dot{\bar{\xi}}\\
\dot{\bar{\eta}}
\end{array}\right)=\frac{1}{\bar{\eta}+\bar{\xi}}\left(\begin{array}{cc}
\mbox{\ensuremath{\mathcal{J}}}_{1} & 1\\
1 & \mathcal{J}_{2}
\end{array}\right)\left(\begin{array}{c}
\bar{\xi}\\
\bar{\eta}
\end{array}\right).
\]
Note that given condition \eqref{eq:Invisible} the eigenvalues of
the matrix in \eqref{eq:SlidNormForm} are real, and given by
\[
\lambda_{1,2}=\frac{1}{2}\left(\mbox{\ensuremath{\mathcal{J}}}_{1}+\mbox{\ensuremath{\mathcal{J}}}_{2}\pm\sqrt{\left(\mbox{\ensuremath{\mathcal{J}}}_{1}-\mbox{\ensuremath{\mathcal{J}}}_{2}\right)^{2}+4}\right),
\]
with associated eigenvectors
\[
v_{1,2}=\left(\begin{array}{c}
\lambda_{1,2}-\mathcal{J}_{2}\\
1
\end{array}\right).
\]
The first eigenvector points into the sliding region, while the
second eigenvector points into the crossing region. One of the eigenvalues
become zero if $1=\mbox{\ensuremath{\mathcal{J}}}_{1}\mbox{\ensuremath{\mathcal{J}}}_{2}$.
The eigenvalues cannot be equal, but they can have the same magnitude
with different sign if $\mbox{\ensuremath{\mathcal{J}}}_{2}+\mbox{\ensuremath{\mathcal{J}}}_{1}=0$.
Depending on the values of $\mbox{\ensuremath{\mathcal{J}}}_{1},\mbox{\ensuremath{\mathcal{J}}}_{2}$
we have three cases \cite{jc09,cj10}:
\begin{enumerate}
\item If $\mbox{\ensuremath{\mathcal{J}}}_{1},\mbox{\ensuremath{\mathcal{J}}}_{2}<0$
and $\mbox{\ensuremath{\mathcal{J}}}_{1}\mbox{\ensuremath{\mathcal{J}}}_{2}>1$
both eigenvalues are negative, so the singularity attracts trajectories from the sliding region, and repels them into the escaping region. The eigenvector in the sliding
region is weak stable compared to the eigenvector in the crossing
region, which implies that trajectories from the sliding region all approach the weak eigenvector direction, and hence all flow into the singularity.
\item If $\mbox{\ensuremath{\mathcal{J}}}_{1},\mbox{\ensuremath{\mathcal{J}}}_{2}>0$
and $\mbox{\ensuremath{\mathcal{J}}}_{1}\mbox{\ensuremath{\mathcal{J}}}_{2}>1$
both eigenvalues are positive, the singularity is repelling from the
sliding region and attracting from the escaping region.
\item When $\mbox{\ensuremath{\mathcal{J}}}_{1}\mbox{\ensuremath{\mathcal{J}}}_{2}>1$
one eigenvalue is positive, the other negative, the singularity is repelling from the
sliding region and attracting from the escaping region.
\end{enumerate}
In terms of our mechanical system case 1 is the most interesting, because
all sliding trajectories in the vicinity of the singularity will reach
$x^*$. This is illustrated in Figure \ref{fig:2fold}.

The local dynamics is illustrated by a few typical trajectories in
Figure \ref{fig:2fold}. Trajectories repeatedly wind around
the singularity through sequences of switchings and stickings, such that even
the local dynamics can be somewhat complex, and a fuller description
can be found in \cite{jc09,cj10}. For the mechanical
system of interest, simulations show that trajectories spend only
a short amount of time near the singularity. In doing so they can
travel through the singularity itself, and what happens at that point
is non-deterministic: a continuous set of future trajectories are
possible.

The problem of unpredictability could be neglected if it affected
only a single point. In this system that is not the case for two reasons.
Firstly, 
whole sets of solutions evolve into the singularity. Secondly,
even after solutions have left the singularity, the global dynamics
can cause their trajectories to return to the singularity again. This
creates a set of closed solutions that return repeatedly to the singularity,
even though which trajectory they follow within that set cannot be
determined. The motion therefore inhabits a chaotic set generated by non-determinism
at the singularity. A simplification of this feedback mechanism is
sketched in Figure \ref{fig:2fold}(right). In the next section we simulate
this non-deterministic chaos in the mechanical model and give parameters
for which it occurs.


\section{Simulations of the mechanical model}\label{sec:sim}

\begin{figure}
\begin{centering}
\includegraphics[width=\linewidth]{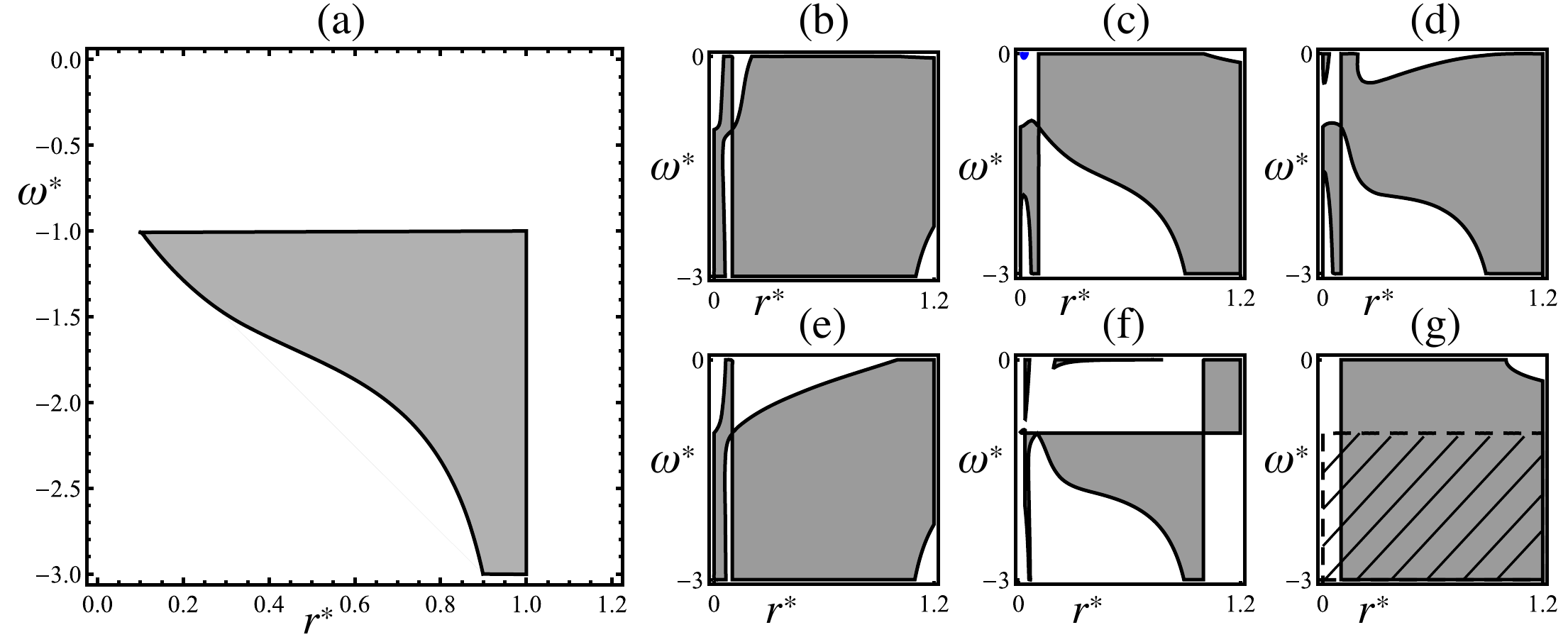}
\caption{\sf The shaded region in panel (a) represents parameter values where the system has
a Teixeira singularity (case 1 in section \ref{sec:teixeira}). This is the overlap of several conditions: (b) $\mathcal{K}^{--}>0$, (c) $\mathcal{K}^{++}<0$, (d) $\mathcal{J}_{1}<0$, (e) $\mathcal{J}_{2}<0$,
(f) $\mathcal{J}_{1}\mathcal{J}_{2}-1>0$, and in (g) showing both $k_{2}>0$ (shaded) and $\kappa>0$ (hatched).}\label{fig:TexRegion}\end{centering}
\end{figure}

A number of conditions must be satisfied for the system (\ref{eq:VF}) to lead to non-determinism and chaos. The presence of a Teixeira singularity exhibiting non-determinism is guaranteed from \eqref{eq:Invisible} if $\mathcal{J}_{1}<0$,
$\mathcal{J}_{2}<0$ and $\mbox{\ensuremath{\mathcal{J}}}_{1}\mbox{\ensuremath{\mathcal{J}}}_{2}>1$ all hold. If we fix the constants $d=1$, $m=1$, $c_{1}=c_{2}=10^{-3}$, $\beta=\frac{1}{20}$,
$r_{0}=10^{-1}$, $\omega_{0}=-1$, $\mu=1$ and $\gamma=-\frac{3}{4}\pi$, then there exist regions of $r^*$ and $\omega^*$ values satisfying these conditions, as illustrated in Figure
\ref{fig:TexRegion}. In general there are large regions in parameter space where the required conditions are satisfied, and while not all will also yield the recurrent dynamics required for chaos, it is not difficult to find parameters that do.


\begin{figure}[h!]
\begin{centering}
\includegraphics[width=0.9\linewidth]{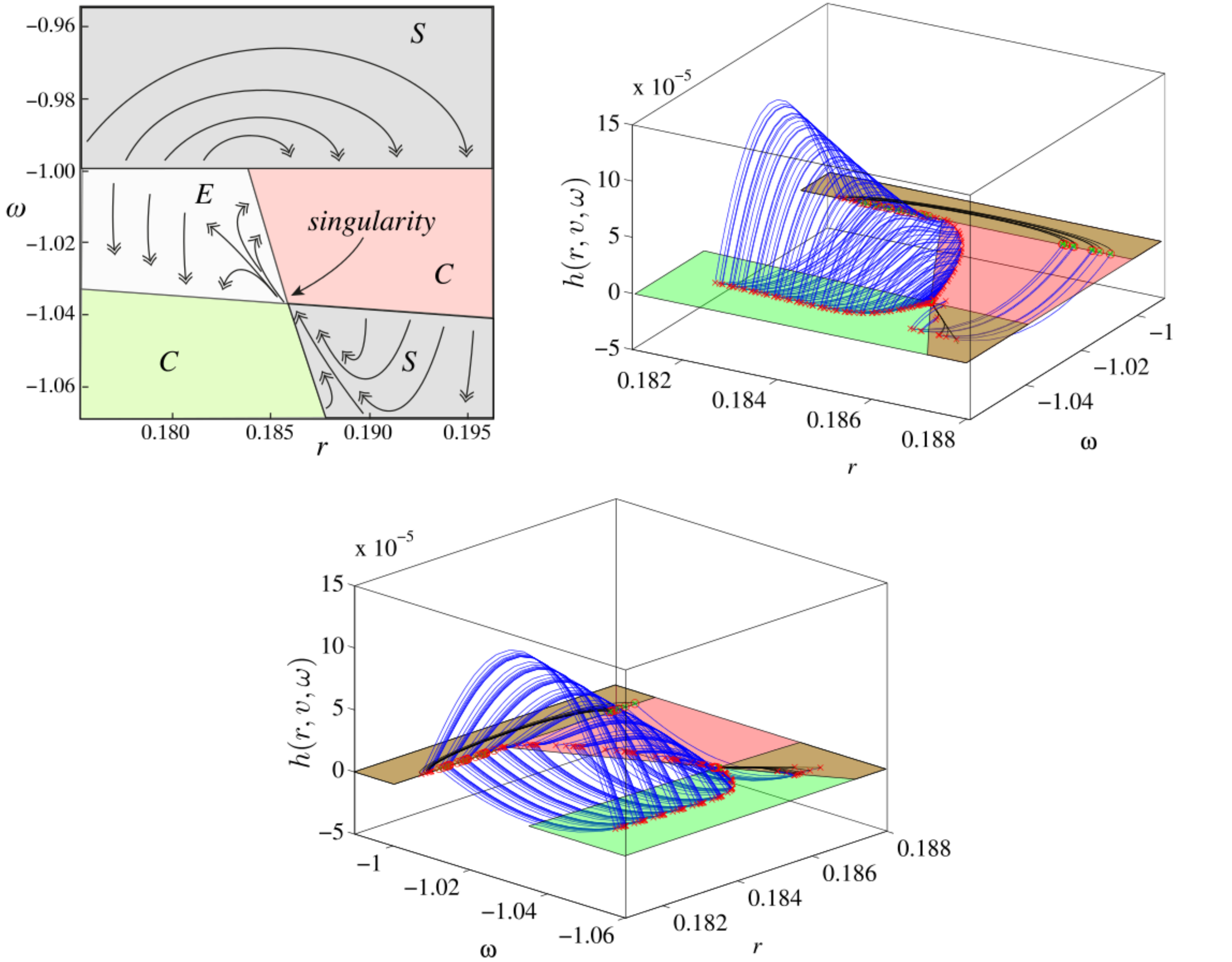}
\vspace{-0.3cm}
\caption{\sf Non-deterministic chaos in the wheel-and-disc assembly for parameter values given in the text. Top left panel: a sketch of the sticking flow, showing sliding regions (S), escaping region (E), and crossing regions (C). A Teixeira singularity is seen at $(r,\omega)=(0.1859,-1.037)$, with the sticking flow passing through it along a weak eigenvector. The 3D blocks show a simulation of different trajectories with initial conditions near the escaping region (regions on $v=0$ correspond to those in the top-left panel; colour online). These represent the different trajectories followed after passing through the singularity. The trajectories wind around the escaping region causing the dense shell seen in $r<0.1859$.
}\label{fig:simula}\end{centering}
\end{figure}

Within the shaded area in Figure \ref{fig:TexRegion}(a) we choose $r^*=0.1859$
and $\omega^*=-1.037$ and use the numerical method described
in \cite{pk08} to compute trajectories of the piecewise-smooth system.
The top left panel of Figure \ref{fig:simula} shows the switching surface
around the singularity on the sliding surface. The sliding surface
is indicated by the grey region where the green and red shaded regions
overlap. The non-overlapping green and red regions are where trajectories
cross the switching surface $\Sigma$. The white area indicates the
escaping surface. It can be seen that trajectories from
the bottom right sliding surface will pass through the singularity,
while the trajectories of the top sliding surface connect the escaping
surface with the middle crossing region.

The full simulation is shown in the top left and bottom panels of Figure \ref{fig:simula}(b).
The figures shows that sticking
trajectories that cross the singularity arrive in the escaping region
and get repelled by it. After leaving the escaping region trajectories
go into a outward spiraling motion and cross the switching surface
until that ceases to exists and arrive at the other sliding surface.
Eventually the solution starts to slip again to end up in the sliding
surface of the singularity to start the cycle again. Because of the
non-determinism at the singularity the cycles are not identical, hence
the attractor is not a periodic orbit, but a full surface of solutions.
The basin of attraction of this set of solutions is not the entire
phase space. In fact it is even possible that the recurrent non-deterministic
phenomena is only transient, because the non-deterministic effect
can also push the trajectory out towards a different attractor.

These trajectories are not simulated actually passing through the singularity,
as such a simulation is impossible by the very definition of the singularity. When a trajectory enters the singularity (as every trajectory in figure \ref{fig:simula} does), the simulation must be stopped and a decision must be made on which of the infinity of possible trajectories to follow. Since all of these pass through the escaping region at least briefly, we simply select a new starting condition near the escaping surface, and by following these we begin to trace out the cone-like form of the non-deterministic chaotic set, and show at least that all such trajectories find their way eventually back into the singularity via the sliding region.

\section{Special cases: singularity without non-determinism}\label{sec:special}

\subsection{The case $\gamma=0$}

%

When $\gamma=0$ the wheel is mounted so that its axis is perpendicular
to the slider. As a result, when $r=0$ the wheel rolls without slipping
around the lower disc in a circle of radius $d$. If $r$ is nonzero
the wheel will typically slip except when the upper and lower discs
rotate in unison, which is when $\omega=\omega_{0}$. Thus the sticking
surface is $h=r(\omega_0-\omega)$, as is found by substituting
$\gamma=0$ in (\ref{h}).

The equations for a two-fold singularity (\ref{2fold}) then
have only one valid solution where $r^*=\kappa$, $\omega^*=\omega_{0}$, and $v^*=\left(d[\frac{k_{2}}{m}(\kappa+r_{0})-\kappa\omega_{0}^{2}]-\frac{c_{1}}{m}\omega_{0}\right)/\left(d\frac{c_{2}}{m}-2\kappa\omega_{0}\right)$.
Although this is a two-fold singularity, it is a special case in which $f(x^*,+\mu)=f(x^*,-\mu)$, since substituting $\gamma=0$ and
$r=-\kappa$ into (\ref{ds1})-(\ref{ds3}) causes the ${\rm sign}(h)$ terms
to drop out. Hence the
system is continuous, and uniquely defined, precisely at the singularity.
This is a special case of the singularity presented in the previous
section, where non-determinism vanishes as $\gamma$ is taken to $0$.

\subsection{The case $\gamma=\pi/2$}


When $\gamma=\pi/2$ the wheel is mounted with its axis lying along
the slider. If the wheel sticks then its speed in the slider $v$
must be equal to the mismatch in the speed of the discs, $d(\omega-\omega_{0})$,
and substituting $\gamma=\pi/2$ into (\ref{h}) gives $h=d(\omega-\omega_{0})-v$.
There are no two-fold singularities because the sticking boundaries form non-intersecting curves.
This is a common situation in simple mechanical models, where the
expressions for $f(x,\pm\mu)$ differ only by a constant $2\mu$,
so the boundaries between sticking and
slipping given by $h=f(x,\pm\mu)=0$ are displaced by an amount $2\mu$. The boundaries can only meet when $\mu=0$, in which case there is no friction force and no discontinuity.

\section{Closing remarks}\label{sec:conc}

Unlike an abstract model in \cite{JeffreyPRL} contrived to demonstrate
the potential for mechanical models exhibiting non-deterministic chaos,
the model presented here is based on fundamental
mechanical principles which themselves have far wider applications.
The essential characteristics of the model are a coupling
of linear and angular motions, and the sudden change of force that
occurs in dry friction. It is therefore reasonable to suggest that
the two-fold singularity may be common wherever such features occur
in physical applications. Uncertainty in the dynamics from a
particular configuration may now be recognized as a local breakdown
of determinism (whether it repeats to give chaotic dynamics or not),
inherent in the physical geometry.

In this local aspect, the non-deterministic chaotic dynamics here
is very different to standard explanations for the phenomenon of wheel
shimmy or flutter, a problem familiar to anyone who has wrestled with
an unwieldy shopping trolley. Current explanations for this tend to
involve smooth deterministic oscillations around the ideal alignment
of a wheel. These do share some aspects with the non-deterministic
effect in this paper, in that violent oscillations of the wheel are
a result of coupling of linear and rotational motion due to nonlinearities
from the wheel (or tyre) elasticity. Where they differ is that our
model involves unpredictability in the transition through a singular
configuration, and it is discontinuities in the dry-friction model,
more severe than nonlinearity, that create such determinism destroying
singularities. (We should note here that our system is not intended
to model shimmy, since it involves a wheel at a fixed angle moving
around a disc, while shimmy involves a wheel with varying angle moving
along a line; our model is designed only to investigate the phenomenon
of non-deterministic chaos).

In the simulations in section \ref{sec:sim}, it is not actually possible to follow the trajectories through the singularity, because by definition their evolution becomes non-deterministic there. Many possibilities can be argued for how to follow the simulation through the singularity, perhaps by smoothing out the sign functions in (\ref{eq:force})-(\ref{eq:moment}), or by choosing one of the possible trajectories probabilistically assuming some stochastic distribution. The smoothing approach produces a stiff model that is very sensitive to the functional form of the smoothing, to choices of parameters and to initial conditions, and sensitive to machine precision used in the simulation. Such simulations were made on a simplified version of this model in \cite{j13capetown}. The stochastic approach has not been studied to date, as there is no obvious way to choose an appropriate distribution based on likely noise input into the system, however such an investigation would be of great interest. Preliminary investigations such as \cite{j13capetown} suggest that, as expected, however noisy or unpredictable the output of such smoothed or stochastic systems, they do indeed constitute a sample of the trajectories observed in the discontinuous model, i.e. in figure \ref{fig:simula}.

\section*{Appendix A. Equations of motion}

\begin{figure}[h!]
\begin{centering}
\includegraphics[width=0.5\textwidth]{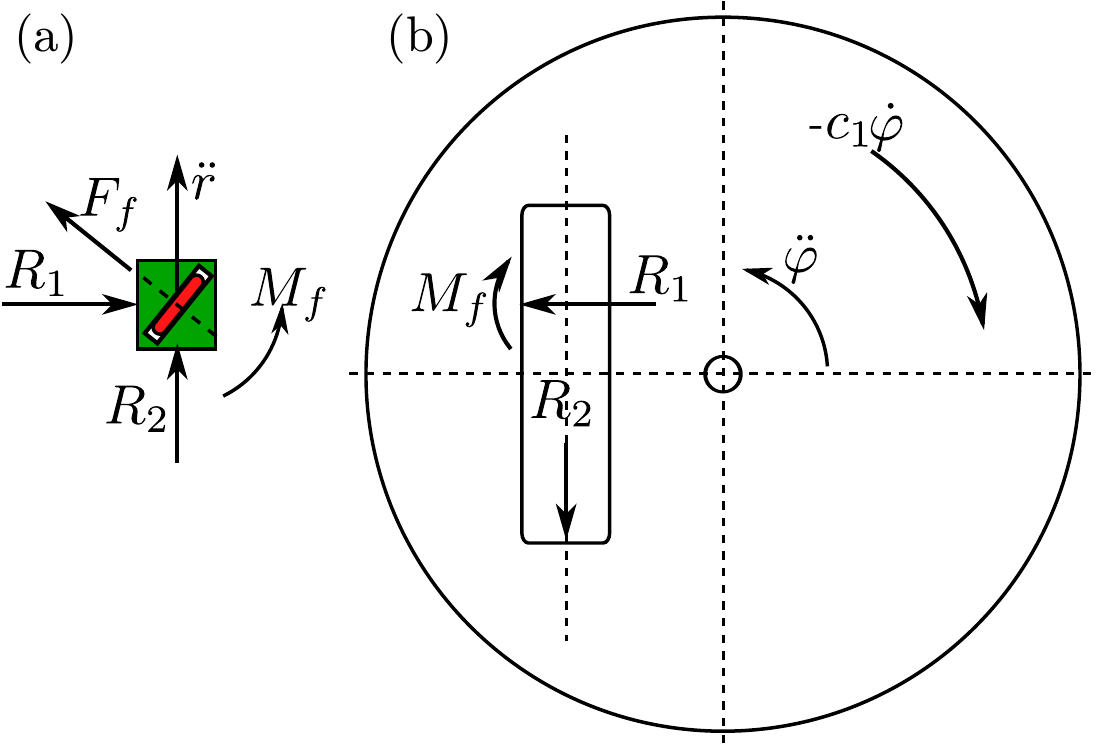}\vspace{-0.11cm}
\caption{\sf Free-body diagrams of the mechanical system.}\label{fig:FBD}\end{centering}
\end{figure}
Consider the system as illustrated by the schematic in Fig.\ \ref{fig:mechmod}.
In order to derive the equations of motions we consider the free-body
diagrams of each part in Fig.\ \ref{fig:FBD}. The generalized coordinates
are the rotation angle $\varphi$ of the top disc and the relative
displacement $r$ of the slider within the disc.
Using these coordinates we calculate the acceleration of the slider in two directions and apply Newton's law to find the equations of motion.
The acceleration of the slider perpendicular to the slit appears in
\begin{equation}
m\left(-r\ddot{\varphi}+d\dot{\varphi}^{2}-2\dot{r}\dot{\varphi}\right)=R_{1} + F(h,g)\cos \gamma\,, \label{eom1}
\end{equation}
where $R_1$ is the reaction force between the slider and the slit and $F(h,g)$ is the friction force between the wheel and the turntable.
Parallel to the slit the motion is governed by
\begin{equation}
m\left(\ddot{r}-d\ddot{\varphi}-r\dot{\varphi^{2}}\right)=R_{2} - F(h,g) \cos \gamma\,, \label{eom2}
\end{equation}
where $R_{2}=-k_{2}(r-r_{0})-c_{2}\dot{r}$ is the reaction force in the spring and damper.
The angular acceleration of the disc is
\begin{equation}
\varTheta\ddot{\varphi}=-c_{1}\dot{\varphi}+dR_{2} + r R_{1} + M(h,g)\,, \label{eom3}
\end{equation}
where $M(h,g)$ is te re-aligning moment of the wheel due to friction. The functions $h$ and $g$
are defined by equations (\ref{h}) and (\ref{g}) with $\omega = \dot{\varphi}$ and $v = \dot{r}$.

Eliminating the reaction forces and rearranging equations (\ref{eom1},\ref{eom2},\ref{eom3}) we arrive at (\ref{ds1},\ref{ds2},\ref{ds3}).
The friction force $F(h,g)$ and moment $M(h,g)$ are derived in the following appendix.

\section*{Appendix B. Friction force on the rolling and slipping wheel}

\begin{figure}
\begin{centering}
\includegraphics{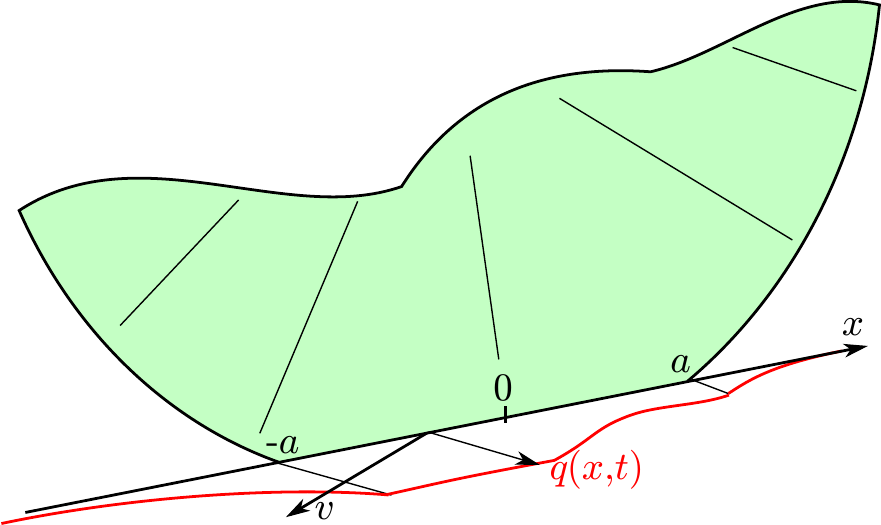}
\caption{\sf Contact forces on a rolling wheel. The ground surface is moving with
relative velocity $v$ under the wheel, whose lateral deformation is described
by $q(x,t)$ within the contact patch spreading from $x=-a$ to $x=a$.}\label{fig:rollingwheel}\end{centering}
\end{figure}

The friction forces acting on the wheel are found as follows. We
use the stretched string tyre model described in \cite{stepan09}
to derive lateral friction forces and the aligning torque. The schematic
of the model can be seen in Figure \ref{fig:rollingwheel}.

The deformation of the tyre is characterized by the lateral movement
of its central line with respect to its resting position by $q(x)$.
We assume that the lateral contact pressure distribution acting on
the tyre is proportional to the displacement that is $p(x)=kq(x)$,
hence the lateral force is
\[
F=k\int_{-a}^{a}q(x)\mathrm{d}x,
\]
where $k$ is the lateral stiffness of the tyre. There is also a torque
that is trying to rotate the wheel back into a position parallel to
the velocity $v$ and that can be calculated for the origin as
\[
M=-k\int_{-a}^{a}xq(x)\mathrm{d}x.
\]
To determine the shape of $q(x)$ we recall that in case of rolling
the contact point has the same velocity as the ground surface. Also,
we assume that the direction of the wheel and ground velocity $v$
has a $\varphi$ angle between them. We are interested in the static
deformation, therefore the following constraint holds
\begin{equation}
0=v\sin\varphi+v\frac{\partial}{\partial x}q(x)\cos\varphi.\label{eq:rollconstr}
\end{equation}
The boundary condition is specified at the leading edge $\sigma\frac{\partial}{\partial x}q(x)=-q(a)$,
where $\sigma$ is a relaxation constant that specifies how the tyre
deformation relaxes when it is not in contact with the ground. Solving
\eqref{eq:rollconstr} we find that $q(x)=q_{0}-x\tan\varphi.$ From
the boundary condition we calculate that $q_{0}=(a+\sigma)\tan\varphi$,
therefore the steady state deformation when the wheel is rolling becomes
\[
q(x)=(a+\sigma-x)\tan\varphi.
\]
To account for a possible slip the distributed friction force $p(x)$
must be limited. We assume static and dynamic friction, where the
maximal static lateral pressure is $\delta$ and the dynamic pressure
is $\rho\delta$, $0<\rho\le1$. This implies that the deformation
of the tyre is limited by the maximal friction pressure $\delta$.
Solving $k\left|q(x_s)\right|=\delta$ for $x_s$, we find that $x_{s}=a+\sigma-\frac{\cot\varphi}{k}$.
Depending on the value of $x_{s}$ three cases are possible:
\begin{enumerate}
\item $x_{s}\le-a$, that is, there is no slip. In this case $F=-2ak(a+\delta)\tan\varphi$, $M=\frac{2}{3}a^{3}k\tan\varphi$
\item $-a<x_{s}<a$ implies partial slip, therefore $F=\delta(2a+\sigma)+\frac{\delta^{2}}{2k}(1-2\rho)\cot\varphi-\frac{1}{2}k\sigma^{2}\tan\varphi$,
\[
M=\frac{k^{3}\sigma^{2}(3a+\sigma)\tan\varphi+\delta^{2}\cot(\varphi)(2\delta\cot\varphi-3k(a+\sigma))}{6k^{2}}-\frac{1}{2}\delta\rho\left(\left(a-\frac{\delta\cot\varphi}{k}+\sigma\right)^{2}-a^{2}\right)
\]
\item $a\le x_{s}$ means complete slip, so that $F=2a\delta$, $M_{z}=0$.
\end{enumerate}
As a simplification we choose $\sigma=0$, $\rho=2/3$, $\delta=(3\kappa)^{-1}$
and $a=3\kappa$, where $\kappa$ is a length scale. This choice of parameters implies that the tyre is infinitely flexible
for bending, its stiffness is the same as the maximum static friction
pressure and $a$ and $\delta$ are chosen such that $F=1$ and $M=\kappa$
when $x_{s}=-a$. Using these values we get
\[
F=\begin{cases}
18k\kappa^{2}\tan\varphi\\
\frac{4}{3}-\frac{\cot\varphi}{54k\kappa^{2}}
\end{cases},\; M=\begin{cases}
18k\kappa^{3}\tan\varphi & \mbox{if}\;0\le\varphi<\cot^{-1}\left(18k\kappa^{2}\right)\\
\frac{\cot\varphi}{18k\kappa} & \mbox{if}\;\cot^{-1}\left(18k\kappa^{2}\right)\le\varphi\le\frac{\pi}{2}
\end{cases}.
\]

\begin{figure}
\begin{centering}
\includegraphics[width=0.49\linewidth]{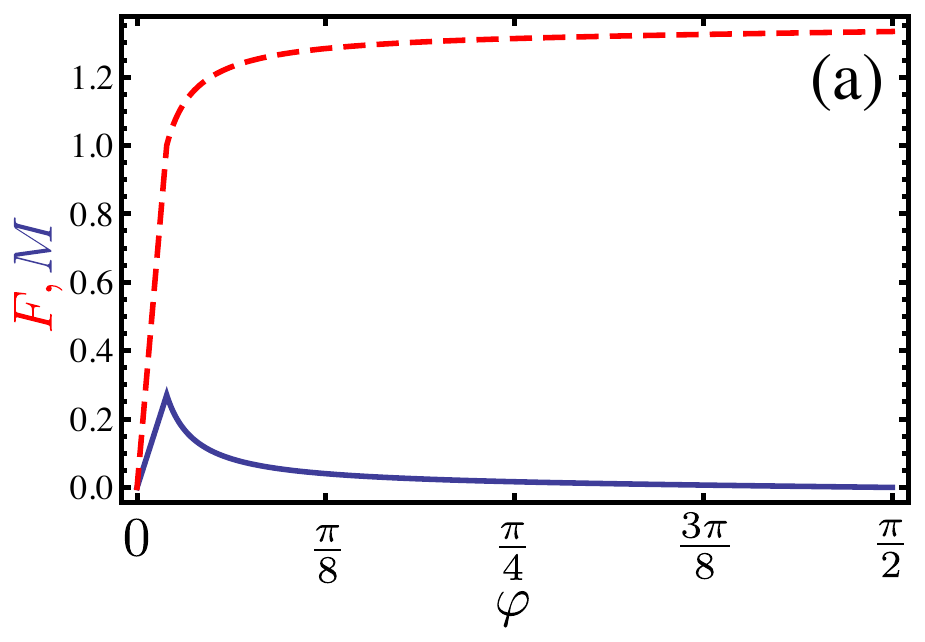}\includegraphics[width=0.49\linewidth]{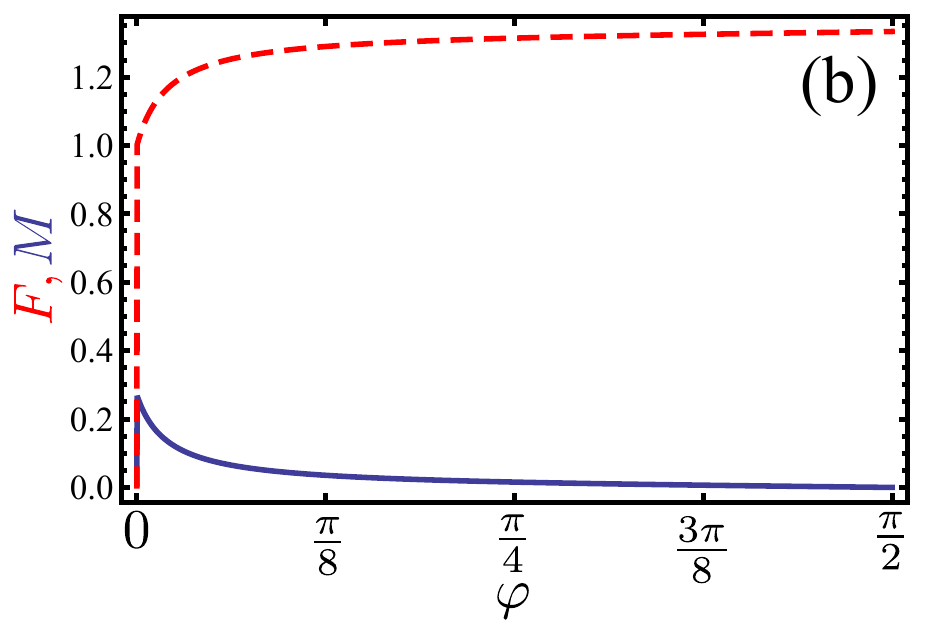}
\caption{\sf Plot of the friction force (dashed red lines) and moment (solid blue lines) as a function
of slip angle $\varphi$. Parameters for the left panel are $k=1$, $\delta=1$, $\sigma=1$
and $a=1$, for right panel are $k=10^{3}$, $\delta=10^{-3}$, $\sigma=10^{-3}$
and $a=1$.}\label{fig:FrictionPlot}\end{centering}
\end{figure}

Choosing some parameters the friction force and torque is plotted
in Fig. \ref{fig:FrictionPlot}(a). Note that this characteristic is
rather similar to Pacejka's empirical magic formula. It can be seen
that the initial part of the $F$ and $M$ curves have a steep gradient.
In order to impart the non-smooth phenomena we replace this initial
section by a discontinuity as shown in Fig. \ref{fig:FrictionPlot}(b).
The scaling is $\varphi\to\cot^{-1}\left(18k\kappa^{2}\right)+\varphi\left(1-\frac{2\cot^{-1}\left(18k\kappa^{2}\right)}{\pi}\right)$.

\bibliography{grazcat}
\bibliographystyle{plain}

\end{document}